\documentclass[11pt]{article} \topmargin=-0.5cm
\usepackage{amssymb}

   \csname @addtoreset\endcsname{equation}{section}
    \csname @addtoreset\endcsname{figure}{section}
    \csname @addtoreset\endcsname{table}{section}

\newcounter{thanksnum}
\def\thanksnumber#1
{\setcounter{thanksnum}{\value{footnote}}\setcounter{footnote}{#1}%
                     \addtocounter{footnote}{-1}\footnotemark
                     \setcounter{footnote}{\value{thanksnum}}}
\def\newtheoremz#1{\@ifnextchar[{\@othmz{#1}}{\@nthmz{#1}}}

\def\@nthmz#1#2{%
\@ifnextchar[{\@xnthmz{#1}{#2}}{\@ynthmz{#1}{#2}}}

\def\@xnthmz#1#2[#3]{\expandafter\@ifdefinable\csname #1\endcsname
{\@definecounter{#1}\@addtoreset{#1}{#3}%
\expandafter\xdef\csname the#1\endcsname{\expandafter\noexpand
  \csname the#3\endcsname \@thmcountersepz \@thmcounterz{#1}}%
\global\@namedef{#1}{\@thmz{#1}{#2}}\global\@namedef{end#1}{\@endtheoremz}}}

\def\@ynthmz#1#2{\expandafter\@ifdefinable\csname #1\endcsname
{\@definecounter{#1}%
\expandafter\xdef\csname the#1\endcsname{\@thmcounterz{#1}}%
\global\@namedef{#1}{\@thm{#1}{#2}}\global\@namedef{end#1}{\@endtheoremz}}}

\def\@othmz#1[#2]#3{\expandafter\@ifdefinable\csname #1\endcsname
  {\global\@namedef{the#1}{\@nameuse{the#2}}%
\global\@namedef{#1}{\@thmz{#2}{#3}}%
\global\@namedef{end#1}{\@endtheoremz}}}

\def\@thmz#1#2{\refstepcounter
    {#1}\@ifnextchar[{\@ythmz{#1}{#2}}{\@xthmz{#1}{#2}}}

\def\@xthmz#1#2{\@begintheoremz{#2}{\csname the#1\endcsname}\ignorespaces}
\def\@ythmz#1#2[#3]{\@opargbegintheoremz{#2}{\csname
       the#1\endcsname}{#3}\ignorespaces}

\def\@thmcounterz#1{\noexpand\arabic{#1}}
\def\@thmcountersepz{.}
\def\@begintheoremz#1#2{ \trivlist \item[\hskip \labelsep{\bf #1\ #2}]}
\def\@opargbegintheoremz#1#2#3{ \trivlist
      \item[\hskip \labelsep{\bf #1\ #2\ (#3)}]}
\def\@endtheoremz{\endtrivlist}

\newtheorem{theorem}{Theorem}[section]
\newtheorem{lemma}{Lemma}[section]

\newtheorem{proposition}{Proposition}[section]
\newtheorem{corollary}{Corollary}[section]
\newtheorem{condition}{Condition}[section]
\newtheorem{definition}{Definition}[section]
\newtheorem{remark}{Remark}[section]

\def\defi{\stackrel{{\scriptscriptstyle \Delta}}{=}}

\def\a{\alpha}
\def\d{\delta}
\def\o{\omega}
\def\O{\Omega}

\def\F{{\cal F}}
\def\w{\widehat}

\def\esssup{\mathop{\rm ess\, sup}}

\def\R{{\bf R}}
\def\E{{\bf E}}
\def\P{{\bf P}}

\def\Z{{\cal Z}}

\def\L{L}

\def\b{\beta}
\def\s{\delta}
\def\g{\gamma}

\def\ww{\widetilde}
\def\X{{\cal X}}

\def\oo{\bar}
\def\s{\sigma}

\def\p{\partial}
\def\G{\Gamma}

\def\A{{\cal A}}
\def\M{{\cal M}}

\def\L{{\cal L}}

\def\TT{{\cal T}}
\newcommand{\be}{\begin{equation}}
\newcommand{\ee}{\end{equation}}
\newcommand{\bd}{\begin{displaymath}}
\newcommand{\ed}{\end{displaymath}}
\newcommand{\ba}{\begin{array}{ll}}
\newcommand{\ea}{\end{array}}
\newcommand{\baa}{\begin{eqnarray}}
\newcommand{\eaa}{\end{eqnarray}}
\newcommand{\baaa}{\begin{eqnarray*}}
\newcommand{\eaaa}{\end{eqnarray*}}   \font\sm=cmr10

\def\ww{\tilde}

\def\u{ u }

\def\QQ{{ Q}}

\def\Q{{\cal Q}}
\def\CC{{\cal C}}
\title{On forward and backward SPDEs with
non-local boundary conditions }
\author{
Nikolai Dokuchaev\\
 {\sm Department of Mathematics \& Statistics, Curtin
University,}\\ {\sm  GPO Box U1987, Perth, 6845 Western Australia} }
\begin{document}
\maketitle
\begin{abstract}
We study linear stochastic partial differential equations of
parabolic type with non-local in time or mixed in time boundary
conditions. The standard Cauchy condition at the terminal time is
replaced by a condition that  mixes the random values of the
solution at different times, including the terminal time, initial
time and continuously distributed times. For the case of backward
equations, this setting covers  almost surely  periodicity.
Uniqueness, solvability and regularity results for the solutions are
obtained. Some possible applications to portfolio selection are
discussed.
\end{abstract}
\section{Introduction}
Stochastic partial
differential equations (SPDEs)  are well studied  in the existing literature
 for the case of Cauchy  boundary conditions at the initial time or at the terminal time.
Forward parabolic SPDEs  are usually considered with a Cauchy
condition at initial time, and backward parabolic SPDEs are usually
considered with a Cauchy condition at terminal time.  A
 backward SPDE cannot be transformed into a forward
  equation by a simple time
 change. Usually, a
  backward SPDE is solvable in the sense that there exists a
  diffusion term being considered as a part of the solution that
  helps to ensure that the solution is adapted to the driving Brownian
motions.  In addition, there are results for the pairs of forward
and backward
 equations with separate  Cauchy conditions at initial time and
 the terminal time respectively. The results for solvability of regularity of
 forward and backward SPDEs can be found in
Al\'os et al (1999), Bally {\it et al} (1994),  Da Prato and Tubaro
(1996), Du and Tang (2012), Gy\"ongy (1998), Krylov (1999),  Ma and
Yong (1997), Maslowski (1995), Pardoux (1993),
 Rozovskii (1990), Walsh (1986),  Yong and Zhou
 (1999), Zhou (1992),
\par
There are also results for SPDEs with boundary conditions that mix
the solution at different times that may  include initial time and
terminal time. This category includes stationary type solutions for
forward SPDEs (see, e.g., Caraballo {\em et al } (2004),
 Chojnowska-Michalik (19987), Chojnowska-Michalik and Goldys
 (1995), Duan {\em et al} (2003), Mattingly (1999),
Mohammed  {\em et al} (2008), Sinai (1996), and the references
there).  Related results were obtained for periodic solutions of
SPDEs
 (Chojnowska-Michalik (1990), Feng and Zhao (2012), Kl\"unger  (2001)).
 As was mentioned in Feng and Zhao (2012), it is difficult to expect that, in general, a SPDE has a periodic
 in time solution $u(\cdot,t)|_{t\in[0,T]}$ in a usual sense of exact equality  $u(\cdot,t)=u(\cdot,T)$ that holds almost surely.
\par
 The  periodicity  of the solutions of stochastic equations  was usually considered
 in the sense of the distributions. In Feng and Zhao (2012),
 the periodicity was established  in a  stronger sense as a "random
 periodic solution" (see Definition 1.1 from Feng and Zhao (2012)); this definition does not assume the equality  $u(\cdot,t)=u(\cdot,T)$.
For ordinary stochastic equation, some periodic solutions were
studied in Rodkina (1992).

The present paper addresses these and related problems again. We
consider  Dirichlet  condition at the boundary of the state domain;
the equations are of  a parabolic type  and are not necessary
self-adjoint. The standard boundary value Cauchy condition at the
one fixed time is replaces by a non-local in time condition that
mixes in one equation the values of the solution at different times
over given time interval, including the terminal time and
continuously distributed times. This is a novel setting comparing
with the periodic conditions for the distributions, or with
conditions from Kl\"unger (2001) and Feng and Zhao (2012), or with
conditions for expectations from Dokuchaev (2008). These conditions
include, for instance, conditions $\kappa u(\cdot,T)=u(\cdot,0)+\xi$
a.e. with $\kappa\in[-1,1]$ and some given $\xi$ (Theorem
\ref{Th5}). We present sufficient conditions for existence and
regularity of the solutions in $L_2$-setting for forward and
backward SPDEs (Theorems \ref{FWTh1}-\ref{FWTh5} and
\ref{Th3}-\ref{Th5}).
\par
Related existence and regularity problems for backward SPDEs were considered in
 Dokuchaev (2012b). In Dokuchaev (2012c), related backward SPDEs were studied in the framework of 
 the contraction theorem in $L_\infty$-space.

Some possible applications to portfolio selection  problems are
discussed (Section \ref{SecF} and Theorem \ref{ThF}). 
\section{The problem setting and definitions}
 We
are given a standard  complete probability space $(\O,\F,\P)$ and
a right-continuous filtration $\F_t$ of complete $\s$-algebras of
events, $t\ge 0$. We are given also a $N$-dimensional Wiener
process $w(t)$ with independent components;  it is a Wiener
process with respect to $\F_t$.
\par
Assume that we are given an open domain $D\subset\R^n$ such that
either $D=\R^n$ or $D$ is bounded  with $C^2$-smooth boundary $\p
D$. Let $T>0$ be given, and let $Q\defi D\times [0,T]$. \par
 We will study the boundary value problems for forward equations  in $Q$
\baaa 
\label{Fparab1} &&d_tu=(\A u+ \varphi)\,dt +\sum_{i=1}^N
[B_iu+h_i]\,dw_i(t), \quad t\ge 0,
\\\label{Fparab10}
&& u(x,t,\o)\,|_{x\in \p D}=0
\\ &&u(\cdot,0)-\G u=\xi.
\label{Fparab2}
\eaaa
and boundary value problems for backward equations  in $Q$
\baaa 
\label{Bparab1} &&d_tu+(\A u+ \varphi)\,dt +\sum_{i=1}^N
B_i\chi_idt=\sum_{i=1}^N\chi_i(t)dw_i(t), \quad t\ge 0,
\\\label{Bparab10}
&& u(x,t,\o)\,|_{x\in \p D}=0
\\ &&u(\cdot,T)-\G u=\xi.
\label{Bparab2}
\eaaa Here $u=u(x,t,\o)$, $\chi_i=\chi_i(x,t,\o)$,
$h_i=h_i(x,t,\o)$, $\varphi=\varphi(x,t,\o)$, $\xi=\xi(x,\o)$,
 $(x,t)\in Q$,   $\o\in\O$.\par
In these boundary problems, $\G$ is a linear operator that maps
functions defined on $Q\times \O$  to functions defines on $D\times
\O$. The operator $\A$ is defined as  \be\label{A}\A v\defi
 \sum_{i=1}^n\frac{\p }{\p x_i}
\sum_{j=1}^n \Bigl(b_{ij}(x,t,\o)\frac{\p v}{\p x_j}(x)\Bigr)
   +\sum_{i=1}^n f_{i}(x,t,\o)\frac{\p v}{\p x_i }(x)
+\,\lambda(x,t,\o)v(x), \ee where $b_{ij}, f_i, x_i$ are the
components of $b$,$f$, and $x$ respectively, and \be\label{B}
B_iv\defi\frac{dv}{dx}\,(x)\,\beta_i(x,t,\o) +\oo
\beta_i(x,t,\o)\,v(x),\quad i=1,\ldots ,N. \ee
\par
We assume that the functions $b(x,t,\o):
\R^n\times[0,T]\times\O\to\R^{n\times n}$, $\b_j(x,t,\o):
\R^n\times[0,T]\times\O\to\R^n$, $\oo\b_i(x,t,\o):$
$\R^n\times[0,T]\times\O\to\R$, $f(x,t,\o):
\R^n\times[0,T]\times\O\to\R^n$, $\lambda(x,t,\o):
\R^n\times[0,T]\times\O\to\R$,   $\chi_i(x,t,\o): \R^n\times
[0,T]\times\O\to\R$, and $\varphi (x,t,\o): \R^n\times
[0,T]\times\O\to\R$ are progressively measurable with respect to
$\F_t$ for all $x\in\R^n$, and the function $\xi(x,\o):
\R^n\times\O\to\R$ is $\F_0$-measurable for all $x\in\R^n$. In fact,
we will also consider $\varphi$ and $\xi$ from wider classes. In
particular, we will consider generalized functions $\varphi$.
\par
If the functions $b$, $f$, $\lambda$, $\varphi$, $\G$, and $\xi$,
 are deterministic, then  $\chi_i\equiv 0$ and equation (\ref{parab1}) is
deterministic.
\subsection*{Spaces and classes of functions} 
We denote by $\|\cdot\|_{ X}$ the norm in a linear normed space
$X$, and
 $(\cdot, \cdot )_{ X}$ denote  the scalar product in  a Hilbert space $
X$.
\par
We introduce some spaces of real valued functions.
\par
 Let $G\subset \R^k$ be an open
domain, then ${W_q^m}(G)$ denote  the Sobolev  space of functions
that belong to $L_q(G)$ together with the distributional
derivatives up to the $m$th order, $q\ge 1$.
\par
 We denote  by $|\cdot|$ the Euclidean norm in $\R^k$, and $\bar G$ denote
the closure of a region $G\subset\R^k$.
\par Let $H^0\defi L_2(D)$,
and let $H^1\defi \stackrel{\scriptscriptstyle 0}{W_2^1}(D)$ be the
closure in the ${W}_2^1(D)$-norm of the set of all smooth functions
$u:D\to\R$ such that  $u|_{\p D}\equiv 0$. Let $H^2=W^2_2(D)\cap
H^1$ be the space equipped with the norm of $W_2^2(D)$. The spaces
$H^k$ and $W_2^k(D)$ are called  Sobolev spaces, they are Hilbert
spaces, and $H^k$ is a closed subspace of $W_2^k(D)$, $k=1,2$.
\par
 Let $H^{-1}$ be the dual space to $H^{1}$, with the
norm $\| \,\cdot\,\| _{H^{-1}}$ such that if $u \in H^{0}$ then
$\| u\|_{ H^{-1}}$ is the supremum of $(u,v)_{H^0}$ over all $v
\in H^1$ such that $\| v\|_{H^1} \le 1 $. $H^{-1}$ is a Hilbert
space.
\par
Let $C_0(\oo D)$ be the Banach space of all functions $u\in C(\oo
D)$  such that $u|_{\p D}\equiv 0$ equipped with the norm from
$C(\oo D)$.
\par We shall write $(u,v)_{H^0}$ for $u\in H^{-1}$
and $v\in H^1$, meaning the obvious extension of the bilinear form
from $u\in H^{0}$ and $v\in H^1$.
\par
We denote by $\oo\ell _{k}$ the Lebesgue measure in $\R^k$, and we
denote by $ \oo{{\cal B}}_{k}$ the $\sigma$-algebra of Lebesgue
sets in $\R^k$.
\par
We denote by $\oo{{\cal P}}$  the completion (with respect to the
measure $\oo\ell_1\times\P$) of the $\s$-algebra of subsets of
$[0,T]\times\O$, generated by functions that are progressively
measurable with respect to $\F_t$.
\par
 We  introduce the spaces
 \baaa
 &&X^{k}(s,t)\defi L^{2}\bigl([ s,t ]\times\Omega,
{\oo{\cal P }},\oo\ell_{1}\times\P;  H^{k}\bigr), \quad\\ &&Z^k_t
\defi L^2\bigl(\Omega,{\cal F}_t,\P; H^k\bigr),\\
&&\CC^{k}(s,t)\defi C\left([s,T]; Z^k_T\right), \qquad k=-1,0,1,2, \\&&
\X^k_c= L^{2}\bigl([ 0,T ]\times\O,\,
\oo{{\cal P} },\oo\ell_{1}\times\P;\; C^k(\oo D)\bigr),\quad k\ge 0,
\\&& \Z^k_c\defi
L_2(\O,\F_T,\P;C^k(D)),\quad k\ge 0. \eaaa

  The
spaces $X^k(s,t)$ and $Z_t^k(s,t)$  are Hilbert spaces.
 \par
 In addition, we introduce the spaces $$ Y^{k}(s,t)\defi
X^{k}(s,t)\!\cap \CC^{k-1}(s,t), \quad k=1,2, $$ with the norm $ \|
u\| _{Y^k(s,T)}
\defi \| u\| _{{X}^k(s,t)} +\| u\| _{\CC^{k-1}(s,t)}. $
\par
\index{For a set $S$ and a normed space ${\rm X}$, we denote by
${\rm B}(S,{\bf X})$ the set of all bounded functions $x:S\to{\rm
X}$. For a set $S$ and a Banach space ${\rm X}$, we denote by ${\rm
B}(S,{\bf X})$ the Banach space of bounded functions $x:S\to{\rm X}$
equipped with the norm $\|u\|_{{\rm B}}=\sup_{s\in S}\|x(s)\|_{{\rm
X}}$. We introduce space $\Z=Z_T^0\cup B(D\times\O)$ equipped with
the norm $\|z\|_{\Z}= \|z\|_{Z_T^0}+\|z\|_{ B(D\times\O) }$.}

For brevity, we shall use the notations
 $X^k\defi X^k(0,T)$, $\CC^k\defi \CC^k(0,T)$,
and  $Y^k\defi Y^k(0,T)$.

\begin{proposition} 
\label{propL} Let $\zeta\in X^0$,
 let a sequence  $\{\zeta_k\}_{k=1}^{+\infty}\subset
L^{\infty}([0,T]\times\O, \ell_1\times\P;\,C(D))$ be such that all
$\zeta_k(\cdot,t,\o)$ are progressively measurable with respect to
$\F_t$, and let $\|\zeta-\zeta_k\|_{X^0}\to 0$. Let $t\in [0,T]$ and
$j\in\{1,\ldots, N\}$ be given.
 Then the sequence of the
integrals $\int_0^t\zeta_k(x,s,\o)\,dw_j(s)$ converges in $Z_t^0$ as
$k\to\infty$, and its limit depends on $\zeta$, but does not depend
on $\{\zeta_k\}$.
\end{proposition}
\par
{\it Proof} follows from completeness of  $X^0$ and from the
equality
\begin{eqnarray*}
\E\int_0^t\|\zeta_{k}(\cdot,s,\o)-\zeta_m(\cdot,s,\o)\|_{H^0}^2\,ds
=\int_D\,dx\,\E\left(\int_0^t\big(\zeta_k(x,s,\o)-
\zeta_m(x,s,\o)\big)\,dw_j(s)\right)^2.
\end{eqnarray*}
\begin{definition} 
\rm Let $\zeta\in X^0$, $t\in [0,T]$, $j\in\{1,\ldots, N\}$, then we
define $\int_0^t\zeta(x,s,\o)\,dw_j(s)$ as the limit  in $Z_t^0$ as
$k\to\infty$ of a sequence $\int_0^t\zeta_k(x,s,\o)\,dw_j(s)$, where
the sequence $\{\zeta_k\}$ is such  as in Proposition \ref{propL}.
\end{definition}
\subsection*{Conditions for the coefficients}
 To proceed further, we assume that Conditions
\ref{cond3.1.A}-\ref{condK} remain in force throughout this paper.
 \begin{condition} \label{cond3.1.A} The matrix  $b=b^\top$ is
symmetric  and bounded. In addition, there exists a constant
$\d>0$ such that
\be
 \label{Main1} y^\top  b
(x,t,\o)\,y-\frac{1}{2}\sum_{i=1}^N |y^\top\b_i(x,t,\o)|^2 \ge
\d|y|^2 \quad\forall\, y\in \R^n,\ (x,t)\in  D\times [0,T],\
\o\in\O. \ee
\end{condition}
\begin{condition}\label{cond3.1.B}
The functions  $f(x,t,\o)$, $\lambda (x,t,\o)$, $\b_i(x,t,\o)$,
and $\oo\b_i(x,t,\o)$, are bounded.
\end{condition}
\begin{condition}\label{condK}  There exists an integer
$m\ge 0$, a set $\{t_i\}_{i=1}^m\subset[0,T]$, and linear continuous
operators $\oo\G:L_2(Q)\to H^0$, $\oo\G_i:H^0\to H^0$, $i=1,..,N$,
such that the operators $\oo\G:L_2([0,T];{\cal B}_1,\ell_1,H^1)\to
W_2^1(D)$ and $\oo\G_i:H^1\to W_2^1(D)$ are continuous and
$$ \G u=\E\{\oo\G
u+\sum_{i=1}^m \oo\G_iu(\cdot,t_i)\}. $$
\end{condition}
\par
By Condition \ref{condK}, the mapping $\G: Y^1\to Z_T^0$ is linear
and continuous. This condition covers cases when
 $$
\oo\G u=\int_0^Tk_0(t)u(\cdot,t)dt,\quad
\oo\G_iu(\cdot,t_i)=k_iu(\cdot,t_i),$$ where $k_0(\cdot)\in
L_2(0,T)$ and $k_i\in\R$. It covers also $\G$ such that $$ \oo\G
u=\int_0^Tdt\int_Dk_0(x,y,t)u(y,t)dx,\quad
\oo\G_iu(\cdot,t_i)(x)=\int_Dk_i(x,y)u(y,t_i)dy,$$
 where $k_i(\cdot)$ are some
regular enough kernels.
\par We introduce the set of parameters $$
\ba {\cal P} \defi \biggl( n,\,\, D,\,\, T,\,\, \G,\
\delta,\,\,\,\,\\ \esssup_{x,t,\o,i}\Bigl[| b(x,t,\o)|+
|{f(x,t,\o)}|+|\lambda(x,t,\o)|+|\b_i(x,t,\o)|+|\oo\b_i(x,t,\o)|\Bigr].
\ea $$
\par
Sometimes we shall omit $\o$.
\section{Forward SPDEs}
\subsection{The definition of solution}
 We will study the following boundary value problem in $Q$
\begin{eqnarray} 
\label{FWparab1} &&d_tu=(\A u+ \varphi)\,dt +\sum_{i=1}^N
[B_iu+h_i]\,dw_i(t), \quad t\ge 0,
\\\label{FWparab10}
&& u(x,t,\o)\,|_{x\in \p D}=0
\\ &&u(x,0,\o)-\G u(\cdot)=\xi(x,\o).
\label{FWparab2}
\end{eqnarray}
Here $u=u(x,t,\o)$, $\varphi=\varphi(x,t,\o)$, $h_i=h_i(x,t,\o)$,
 $(x,t)\in Q$,   $\o\in\O$.
 \par
We do not exclude an important special case when the functions $b$,
$f$, $\lambda$, $\varphi$, and $\xi$,
 are deterministic,
 and $h_i\equiv 0$, $B_i\equiv 0$ ($\forall i)$. In this case, equation (\ref{FWparab1}) is deterministic.

\begin{definition} 
\label{FWdefsolltion} \rm Let $u\in Y^1$,  $\varphi\in X^{-1}$, and
$h_i\in X^0$. We say that equations
(\ref{FWparab1})-(\ref{FWparab10}) are satisfied if
\begin{eqnarray}
&&u(\cdot,t,\o)-u(\cdot,r,\o)\nonumber
\\  &&\hphantom{xxx}= \int_r^t\big(\A u(\cdot,s,\o)+
\varphi(\cdot,s,\o)\big)\,ds+ \sum_{i=1}^N
\int_r^t[B_iu(\cdot,s,\o)+h_i(\cdot,s,\o)]\,dw_i(s)
\label{FWintur}
\end{eqnarray}
for all $r,t$ such that $0\le r<t\le T$, and this equality is
satisfied as an equality in $Z_T^{-1}$.
\end{definition}
Note that the condition on $\p D$ is satisfied in the  sense that
$u(\cdot,t,\o)\in H^1$ for a.e. \ $t,\o$. Further, $u\in Y^1$, and
the value of  $u(\cdot,t,\o)$ is uniquely defined in $Z_T^0$ given
$t$, by the definitions of the corresponding spaces. The integrals
with $dw_i$ in (\ref{FWintur}) are defined as elements of $Z_T^0$.
The integral with $ds$ in (\ref{FWintur}) is defined as an element
of $Z_T^{-1}$. In fact, Definition \ref{FWdefsolltion} requires for
(\ref{FWparab1}) that this integral must be equal  to an element of
$Z_T^{0}$ in the sense of equality in $Z_T^{-1}$.
\subsection{Existence and regularity results}
\begin{theorem}
\label{FWTh1} There exist a number $\kappa=\kappa({\cal P})>0$ such
that problem (\ref{FWparab1})-(\ref{FWparab2}) has an unique
solution in the class $Y^1$,  for any $\varphi\in X^{-1}$, $h_i\in
X^0$, $\xi\in Z_0^0$, and any $\G$ such that $ \|\G\|\le \kappa$,
where $\|\G\|$ is the norms of the operator $\G: Y^1\to Z_0^T$. In
addition, \be \label{FW3.3} \| u \|_{Y^1}\le C \left(\| \varphi \|
_{X^{-1}}+\|\xi\|_{Z_0^0} +\sum_{i=1}^N\|h_i \|_{X^0}\right), \ee
where  $ C=C(\kappa,{\cal P})>0$ is a constant that depends only on
$\kappa$ and ${\cal P}$.
\end{theorem}
\par
Starting from now and up to the end of this section, we assume that
Condition \ref{FWcondB} holds.
\begin{condition}\label{FWcondB} \begin{enumerate}
\item The domain $D$ is bounded.
The functions  $b(x,t,\o)$, $f(x,t,\o)$, $\lambda (x,t,\o)$,
$\b_i(x,t,\o)$ and $\oo\b_i(x,t,\o)$ are  differentiable in $x$ for
a.e. $t,\o$, and the corresponding derivatives are bounded.
\item Condition \ref{condK} is satisfied with
$\{t_k\}_{k=1}^m\subset(0,T]$. \end{enumerate}
\end{condition}
It follows from this condition that there exist modifications of
$\b_i$ such that the functions $\b_i(x,t,\o)$ are continuous in $x$
for a.e. $t,\o$. We assume that $\b_i$ are such functions.
\begin{theorem}
\label{FWTh3} Let $\F_0$ be  the $\P$-augmentation of the set
$\{\emptyset,\O\}$. Assume that at least one of the following
conditions is satisfied:
\begin{itemize}
\item[(i)] the function $b$ is non-random, or \item[(ii)]
$\b_i(x,t,\o)=0$ for $x\in \p D$, $i=1,...,N$.
\end{itemize}
\par
Further, assume that  problem (\ref{FWparab1})-(\ref{FWparab2}) with
$\varphi\equiv 0$, $h_i\equiv 0$, $\xi\equiv 0$, does not admit
non-zero solutions in the class $Y^1$. Then   problem
(\ref{FWparab1})-(\ref{FWparab2}) has a unique solution $u$ in the
class $Y^1$ for any $\varphi\in X^{-1}$, $h_i\in X^0$, and $\xi\in
H^0$. In addition, \be \label{FW3.5} \| u \|_{Y^1}\le C \left(\|
\varphi \| _{X^{-1}}+\|\xi\|_{H^0}+ \sum_{i=1}^N\|h_i \|_{X^0}
\right), \ee where  $C>0$ is a constant that
 does not depend on $\varphi, h_i$, and  $\xi$.
\end{theorem}
\begin{theorem}
\label{FWTh5} Let the functions $b,f$ and $\lambda$ be non-random
and such that the operator $\A$ can be represented as $$ \A
v=\sum_{i,j=1}^n\frac{\p^2 }{\p x_i \p x_j}
\left(b_{ij}(x,t)v(x)\right)+\sum_{i=1}^n\frac{\p}{\p x_i }\left( \w
f_i(x,t)v(x)\right)+\w\lambda(x,t)v(x), $$ where $\w\lambda(x,t)\le
0$, and where $\w f_i$ are bounded functions. Further, let $$ \G
u=\E\Bigl\{\int_0^Tk_0(t)u(\cdot,t)dt+\sum_{i=1}^m
k_iu(\cdot,t_i)\Bigl|\F_0\Bigr\}, $$ where $t_i>0$, and where
$k_i\in\R$, $k_0(\cdot)\in L_2(0,T)$ are such that $$
\int_0^T|k_0(t)|dt+\sum_{i=1}^m |k_i|\le 1. $$ Then problem
(\ref{FWparab1})-(\ref{FWparab2}) has a unique solution $u$ in the
class $Y^1$ for any  $\varphi\in X^{-1}$, $h_i\in X^0$, and $\xi\in
H^0$. In addition, (\ref{FW3.5}) holds with a constant $C>0$ that
 does not depend on $\varphi, h_i$, and  $\xi$.
\end{theorem}
\par
The following corollary is a special case of Theorem \ref{FWTh5} for
deterministic parabolic equation with the boundary condition that
covers the condition of periodicity.
\begin{corollary}
\label{FWcor1} Under the assumptions of Theorem \ref{FWTh5}, for any
$k\in[-1,1]$, the deterministic boundary value problem $$ \frac{\p
u}{\p t}=\A u+\varphi,\quad u|_{\p D}=0,\quad u(x,0)-ku(x,T)\equiv
\Phi(x) $$ has a unique solution $u\in C([0,T];H^0)\cap
L_2([0,t],{\cal B}_1,\ell_1,H^1)$ for any $\Phi\in H^0$, $\varphi\in
L_2(Q)$,  and $$ \|u(\cdot,t)\|_{Y^1}\le
C(\|\Phi\|_{H^0}+\|\varphi\|_{L_2(Q)}), $$ where $C>0$ is a constant
that does not depend on $\Phi$  and $\varphi$.
\end{corollary}
The classical result about well-posedness of the Cauchy condition at
initial time corresponds to the special case of $k=0$.
\section{Backward SPDEs}\subsection{The definition of solution} For backward SPDEs, we will study the following boundary
value problem in $Q$
\begin{eqnarray} 
\label{parab1} &&d_tu+(\A u+ \varphi)\,dt +\sum_{i=1}^N
B_i\chi_idt=\sum_{i=1}^N\chi_i(t)dw_i(t), \quad t\ge 0,
\\\label{parab10}
&& u(x,t,\o)\,|_{x\in \p D}=0
\\ &&u(\cdot, T)-\G u(\cdot)=\xi.
\label{parab2}
\end{eqnarray}
Here $u=u(x,t,\o)$, $\varphi=\varphi(x,t,\o)$,
$\chi_i=\chi_i(x,t,\o)$,
 $(x,t)\in Q$,   $\o\in\O$.\par
   In (\ref{parab2}), $\G$ is a linear operator that maps functions
defined on $Q\times \O$  to functions defined on $D\times \O$. For
instance, the case where  $\G u=u(\cdot,0)$ is not excluded; this
case corresponds to the periodic type boundary condition
$u(\cdot,T)-u(\cdot,0)=\xi.$
\begin{definition} 
\label{defsolltion} \rm Let $u\in Y^1$, $\chi_i\in X^0$,
$i=1,...,N$, and $\varphi\in X^{-1}$. We say that equations
(\ref{parab1})-(\ref{parab10}) are satisfied if \baaa
&&u(\cdot,t,\o)=u(\cdot,T,\o)+ \int_t^T\big(\A u(\cdot,s,\o)+
\varphi(\cdot,s,\o)\big)\,ds \ \nonumber
\\&&\hphantom{xxx}+ \sum_{i=1}^N
\int_t^TB_i\chi_i(\cdot,s,\o)ds-\sum_{i=1}^N
\int_t^T\chi_i(\cdot,s)\,dw_i(s)
\label{intur} \eaaa for all $r,t$ such that $0\le r<t\le T$, and
this equality is satisfied as an equality in $Z_T^{-1}$.
\end{definition}
Note that the condition on $\p D$ is satisfied in the  sense that
$u(\cdot,t,\o)\in H^1$ for a.e. \ $t,\o$. Further, $u\in Y^1$, and
the value of  $u(\cdot,t,\o)$ is uniquely defined in $Z_T^0$ given
$t$, by the definitions of the corresponding spaces. The integrals
with $dw_i$ in (\ref{intur}) are defined as elements of $Z_T^0$. The
integral with $ds$ in (\ref{intur}) is defined as an element of
$Z_T^{-1}$. In fact, Definition \ref{defsolltion} requires for
(\ref{parab1}) that this integral must be equal  to an element of
$Z_T^{0}$ in the sense of equality in $Z_T^{-1}$.
\subsection{Existence and regularity results}
Starting from now and up to the end of this section, we assume that
Condition \ref{condB} holds.
\begin{condition}\label{condB}
\begin{enumerate}
\item Condition \ref{FWcondB}(i) is satisfied.
\item $\b_i(x,t,\o)=0$ for $x\in \p D$,
$i=1,...,N$. \item $\F_0$ is the $\P$-augmentation of the set
$\{\emptyset,\O\}$.
\item
Condition \ref{condK} is satisfied with
$\{t_i\}_{i=1}^m\subset[0,T)$.
\end{enumerate}
\end{condition}

In particular, it follows from this condition that there exist
modifications of $\b_i$ such that the functions $\b_i(x,t,\o)$ are
continuous in $x$ for a.e. $t,\o$. We assume that $\b_i$ are such
functions.
\par
Note that the assumptions on $\G$ imposed in Condition \ref{condB}
allows to consider $\G u=u(\cdot,0)$, i.e., the periodic type
boundary conditions $u(\cdot,T)=u(\cdot,0)$ a.s.; it suffices to
assume that  $t_1=0$, $\oo G_1$ is identical operator, $\oo\G=0$,
$\oo G_i=0$, $i>1$.
\begin{theorem}
\label{Th3}   Assume that problem (\ref{parab1})-(\ref{parab2}) with
$\varphi\equiv 0$, $\xi\equiv 0$, does not admit non-zero solutions
$(u,\chi_1,...,\chi_N)$  in the class $Y^1\times
(X^0)^N$. Then   problem (\ref{parab1})-(\ref{parab2}) has a unique
solution $(u,\chi_1,...,\chi_N)$  in the class $Y^1\times (X^0)^N$,
 for any $\varphi\in X^{-1}$, and $\xi\in H^0$. In addition, \be
\label{3.5} \| u \|_{Y^1}+\sum_{i=1}^N\|\chi_i\|_{X^0}\le C \left(\|
\varphi \| _{X^{-1}}+\|\xi\|_{H^0}\right), \ee where $C>0$
 does not depend on $\varphi$ and  $\xi$.
\end{theorem}
\par
Up to the end of this section, we assume that the following condition is satisfied.
\begin{condition}\label{condTh5}
\begin{enumerate}
 \item $\oo\b_i\equiv 0$.
 \item The functions $b,f$ and
$\lambda$ are such that the operator $\A$ can be represented as $$ \A
v=\sum_{i,j=1}^nb_{ij}(x,t,\o)\frac{\p^2 v}{\p x_i \p x_j}(x)
+\sum_{i=1}^n\w f_i(x,t,\o)\frac{\p v}{\p x_i
}(x)+\w\lambda(x,t,\o)v(x),
$$
where the functions  $\w f(x,t,\o)$, $\w \lambda (x,t,\o)$, and $\b_i(x,t,\o)$
are bounded and  are differentiable
in $x$ for a.e. $t,\o$, and the corresponding derivatives are
bounded. \item $b\in \X_c^3$, $\w f\in\X_c^2$,
$\w\lambda\in\X^1_c$, $ \b_i\in\X_c^3$.
\end{enumerate}
\end{condition}
\begin{theorem}
\label{Th5}  Let $\w\lambda(x,t,\o)\le 0$ a.e., and let $\G u=\kappa\u(\cdot,0)$, where
$\kappa\in[-1,1]$, i.e, boundary condition (\ref{parab2}) is \baa
u(\cdot,T)-\kappa u(\cdot,0)=\xi. \label{kappa} \eaa
 Then problem
(\ref{parab1})-(\ref{parab10}),(\ref{kappa}) has
 a unique solution
$(u,\chi_1,...,\chi_N)$ in the class $Y^1\times(X^0)^N$ for any
$\varphi\in X^{-1}$ and $\xi\in  Z^0_T$. In addition, (\ref{3.5}) holds
with $C>0$ that
 does not depend on $\varphi$ and  $\xi$.
\end{theorem}
\par
Conditions (\ref{kappa}) were
introduced in Dokuchaev (1994) for deterministic parabolic equations  (see also Theorem 2.2 from Dokuchaev
(2004)).
\section{Some applications: portfolio selection problems}\label{SecF}
Theorem \ref{Th5} can be applied to portfolio selection  for
 continuous time diffusion market model,
 where the price dynamic is described by Ito stochastic differential equations.
Examples of these models can be found in, e.g.,  \index{\cite{K}}
{Karatzas and   Shreve (1998)}.
\par
We consider the following   model of a securities market consisting
of a risk free bond or bank account with the price $B(t) $, ${t\ge
0}$, and
 a risky stock with the price $S(t)$, ${t\ge 0}$. The prices of the stocks evolve
 as \be \label{S} dS(t)=S(t)\left(a(t)dt+\s(t) dw(t)+\w\s(t)d\w w(t)\right), \quad
t>0, \ee where $(w(t),\w w(t))$ is a Wiener process, $a(t)$ is a
appreciation rate, $(\s(t),\w s(t))$  is a vector of volatility
coefficients. The initial price $S(0)>0$ is a  given deterministic
constant. The price of the bond evolves as \baaa \label{Bond}
B(t)=e^{rt}B(0), \eaaa where $B(0)$ is a  given constant, $r\ge 0$
is a short rate. For simplicity, we assume that $r=0$ and
$B(t)\equiv B(0)$.
\par
 We  assume that $(w(\cdot),\w w(t))$ is a standard Wiener process on  a
given standard probability space $(\Omega,\F,\P)$, where $\Omega$ is
a set of elementary events, $\F$ is a complete $\s$-algebra of
events, and $\P$ is a probability measure.
\par
 Let $\F_t$
be the filtration generated by $w(t)$, and let $\w\F_t$ be the
filtration generated by $(w(t),\w w(t))$. In particular, we assume
that $\F_0$ and $\w \F_0$ are trivial $\s$-algebras, i.e., they are
the $\P$-augmentations of the set $\{\emptyset,\Omega\}$.
\par
We assume that the processes  $a(t)$, $\s(t)$,  $\w\s(t)$,
$\s(t)^{-1}$, and $\w\s(t)^{-1}$ are measurable, bounded and
$\F_t$-adapted.  \index{We assume that the processes $a(t)$,
$\s(t)$, and $\s(t)^{-1}$ are bounded, the process $(S(t),a(t))$ is
$\F_t$-adapted, and  the process $\s(t)$ is continuous and
deterministic. In particular, this means that the process the
process $a(t)$ can be random and that the process $(S(t),a(t))$ is
currently observable.}\subsubsection*{Strategies for
bond-stock-options market} The rules for the operations of the
agents on the market define the class of admissible strategies where
the optimization problems have to be solved.
\par
Let $X(0)>0$ be the initial wealth at time $t=0$ and let $X(t)$ be
the wealth at time $t>0$.
\par
We assume that the wealth $X(t)$ at time $t\in[0,T]$ is
\begin{equation}
\label{X} X(t)=\b(t)B(t)+\g(t)S(t).
\end{equation}
Here $\b(t)$ is the quantity of the bond portfolio, $\g(t)$ is the
quantity of the stock  portfolio, $t\ge 0$. The pair $(\b(\cdot),
\g(\cdot))$ describes the state of the bond-stocks securities
portfolio at time $t$. Each of  these pairs is  called a strategy.
\par  A pair $(\b(\cdot),\g(\cdot))$  is said to be an admissible
strategy if the processes $\b(t)$ and $\g(t)$ are progressively
measurable with respect to the filtration $\w\F_t$.

 In particular, the agents are not supposed to know the
future (i.e., the strategies have to be adapted to the flow of
current market information).

Let $\P_*$ be an equivalent probability  measure such that $S(t)$ is
a martingale under $\P_*$. By the assumptions on $(a,\s,\w s)$, this
measure exists and is unique.

A pair $(\b(\cdot),\g(\cdot))$  is said to be an admissible
self-financing strategy, if
 \baaa
 X(t)=X(0)+\int_0^t\g(s)dS(s),
 \eaaa
 and $\E_*X(T)^2<+\infty$.
\subsection*{A portfolio selection problem}
In portfolio theory, a typical problem is constructing a portfolio
strategy with certain desirable properties. It will be demonstrated
below that Theorem \ref{Th5} can be applied to this problem.

Let us consider the following example.

Let  $s_L\in(0,S(0))$, $s_U\in(S(0),+\infty)$. Let $D=(s_L,s_U)$.
and  $\xi\in L_{\infty}(\O,\F_T,\P,C_0(\oo D))$ be given.

Let us consider  problem (\ref{parab1})-(\ref{parab2}) with $n=N=1$,
$D=(s_L,s_U)$, $\varphi\equiv 0$, and with
 \baaa &&\A
v=\frac{1}{2}(\s(t)^2+\w\s(t)^2)x^2\frac{\p^2 u}{\p x^2}(x), \qquad
B_1v\defi x\s(t)\frac{dv}{dx}\,(x),\qquad (\G u)(x)=u(x,0). \eaaa In
other words, we consider the following problem  \baa
&&d_tu(x,t)+\frac{1}{2}(\s(t)^2+\w\s(t)^2)
x^2\frac{\p^2 u}{\p x^2}(x,t)+\s(t)x\frac{d\chi }{dx}(x,t)=\chi(x,t) dw(t),\quad t<T,\hphantom{xxx}\label{f1}\\
&& u(s_L,t)=u(s_U,t)=0,\label{f2}\\&&
u(x,T)=u(x,0)+\xi(x).\label{f3} \eaa The assumptions of Theorem
\ref{Th5} are satisfied for this problem. By this theorem, there
exists a unique solution $u(x,t,\o): [s_L,s_U]\times[0,T]\times
\O\to\R$ of problem   (\ref{f1})-(\ref{f3})  such that  $u\in Y^1$.

\par
Let $\tau=\inf\{t>0: S(t)\notin D\}$.
\begin{theorem}\label{ThF}
Let  \baa
 X(t,x)=u(S(t\land\tau),t\land\tau),\label{Xmar}\eaa where $x\in D$, $t\in[0,T]$
and where $S(t)$ is defined by (\ref{S}) given that $S(0)=x$. Then
$\E_*X(T,x)^2<+\infty$, and the process $X(t,x)$  represents the
wealth generated by some self-financing strategy given that
$S(0)=x$. In addition, \baa
&X(T,y)=X(0,x)+\xi(x)\quad&\hbox{if}\quad \tau>T,\quad
S(T)y=S(0)x.\label{XT}\eaa
\end{theorem}
\par
The portfolio described in Theorem \ref{ThF} has the following
attractive feature: with a positive $\xi$, it ensures a systematic
gain when $\tau>T$ for the case of stagnated marked prices. The
event $\tau<T$  can be considered as an extreme event if $s_L$ is
sufficiently small and $s_U$ is sufficiently large.
\par
Note that the assumption  that the process $(a(t),\s(t))$ is
$\F_t$-adapted was used to ensure existence of $u$. A more general
model where this process is $\w\F_t$-adapted leads to a degenerate
SPDE in bounded domain where Condition \ref{cond3.1.A} is not
satisfied. This case is not covered  by Theorem \ref{Th5}. An
example  of portfolio selection based on a degenerate backward SPDE
in entire space was considered in Ma and Yong (1997).
\section{Proofs}
{\em Proof of Theorems \ref{FWTh1}-\ref{FWTh5}} can be found in
Dokuchaev (2008).

Let $s\in (0,T]$, $\varphi\in X^{-1}$ and $\Phi\in Z^0_s$. Consider
the problem \be \label{4.1}
\begin{array}{ll}
d_tu+\left( \A u+ \varphi\right)dt +
\sum_{i=1}^NB_i\chi_i(t)dt=\sum_{i=1}^N\chi_i(t)dw_i(t), \quad t\le s,\\
u(x,t,\o)|_{x\in \p D}, \\
 u(x,s,\o)=\Phi(x,\o).
\end{array}
 \ee
 \par
The following lemma represents an analog of the so-called "the first
energy inequality", or "the first fundamental inequality" known for
deterministic parabolic equations (see, e.g., inequality (3.14) from
Ladyzhenskaya (1985), Chapter III).
\begin{lemma}
\label{lemma1} Assume that Conditions \ref{cond3.1.A}--\ref{condK}
are satisfied.  Then problem (\ref{4.1}) has an unique solution a
unique solution $(u,\chi_1,...,\chi_N)$ in the class
$Y^1\times(X^0)^N$  for any $\varphi\in X^{-1}(0,s)$, $\Phi\in
Z_s^0$, and \be \label{4.2} \| u
\|_{Y^1(0,s)}+\sum_{i=1}^N\|\chi_i\|_{X^0}\le C \left(\| \varphi \|
_{X^{-1}(0,s)}+\|\Phi\|_{Z^0_s}\right), \ee where $C=C({\cal P})$
does not depend on $\varphi$ and $\xi$.
\end{lemma}
(See, e.g., Dokuchaev (1991) or Theorem 4.2 from Dokuchaev (2010)).
\par Note that the solution $u=u(\cdot,t)$
is continuous in $t$ in $L_2(\O,\F,\P,H^0)$, since
$Y^1(0,s)=X^{1}(0,s)\!\cap \CC^{0}(0,s)$.
\par
Introduce  operators $L_s:X^{-1}(0,s)\to Y^1(0,s)$ and
$\L_s:Z^0_s\to Y^1(0,s)$, such that $u=L_s\varphi+\L_s\Phi,$ where
 $(u,\chi_1,...,\chi_N)$ is the solution of  problem (\ref{4.1})  in the class
$Y^2\times(X^1)^N$. By Lemma \ref{lemma1}, these linear operators
are continuous.
\par
Introduce   operators $\Q:Z_T^0\to Z_T^0$ and $\TT:X^{-1}\to Z_T^0$
such that $\Q\Phi+\TT\varphi= \G u$,
    where
$u$ is the solution in $Y^1$ of   problem (\ref{4.1}) with $s=T$,
$\varphi\in X^{-1}$, and $\Phi\in Z_T^0$. It is easy to see that
these operators are linear and continuous.
\par
 For brevity, we denote
$u(\cdot,t)=u(x,t,\o)$. Clearly,
 $u\in Y^1$ is the solution of   problem
(\ref{parab1})-(\ref{parab2}), if \baaa
&&u=L_T\varphi+\L_Tu(\cdot,T), \\
&&u(\cdot,T)-\G u=\xi. \eaaa Since $ \G u=\Q u(\cdot,T)+\TT\varphi$,
we have  \baaa
 &&u(\cdot,T)-\Q u(\cdot,T)-\TT\varphi
=\xi. \eaaa Clearly,  $\|\Q\|\le \|\G\|\|\L_T\|$, where $\|\Q\|$,
$\|\G\|$, and  $\|\L_T\|$, are the norms of the operators $\Q:
Z_T^0\to Z_T^0$, $\G: Y^1\to Z_T^0$, and $\L_0: Z_T^0\to Y^1$,
respectively. Since the operator $\Q: Z_T^0\to Z_T^0$ is continuous,
the operator $(I-\Q)^{-1}:Z_T^0\to Z_T^0$ is continuous for small
enough $\|\Q\|$, i.e. for a small enough $\kappa>0$. Hence $$
u(\cdot,T)=(I-\Q)^{-1}(\xi +\TT\varphi), $$ and \baa \label{4.3}
u&=&L_T\varphi+\L_T u(\cdot,T)\nonumber\\&=&L_T\varphi+\L_T
(I-\Q)^{-1}(\xi+\TT\varphi). \eaa
\par
Starting from now, we assume that Condition \ref{condB} is satisfied, in addition
to Conditions \ref{cond3.1.A}-\ref{condK}.
\par
The following lemma represents an analog of the so-called "the
second energy inequality", or "the second fundamental inequality"
known for the deterministic parabolic equations (see, e.g.,
inequality (4.56) from Ladyzhenskaya (1985), Chapter III).
\begin{lemma}\label{Th3.1.1}
Problem (\ref{4.1})   has a unique solution $(u,\chi_1,...,\chi_N)$
in the class $Y^2\times(X^1)^N$ for any $\varphi \in X^0$, $\Phi \in
Z_T^1$, and \be \| u \|_{{Y}^2}+\sum_{i=1}^N\|\chi_i\|_{X^1}
 \le   C\left(  \|\varphi \|_{X^0} +
\|\Phi \|_{Z_T^1} \right), \label{3.1.3} \ee where $C>0$ does not
depend on $\varphi$ and $\Phi$; it depends on ${\cal P}$ an on the
supremums of the derivatives listed in Condition \ref{condB}(ii).
\end{lemma}
\par
The lemma above represents a reformulation of Theorem 3.1. from Du
and Tang (2012), or Theorem 3.4 from Dokuchaev (2010) or Theorem 4.3
from Dokuchaev (2012a). In the cited papers, this result was obtained
under some strengthened version of Condition \ref{cond3.1.A}; this
was restrictive.  In Du and Tang (2012), this result was obtained
without this restriction, i.e., under Condition \ref{cond3.1.A}
only.
\begin{remark}\label{remDu} {  Thanks  to Theorem
3.1. from Du and Tang (2012), Condition 3.5 from Dokuchaev (2011)
and Condition 4.1 from Dokuchaev (2012a)  can be replaced by less
restrictive Condition \ref{cond3.1.A}; all results from Dokuchaev
(2011,2012a) are still valid.}
\end{remark}
\par
\begin{lemma}
\label{lemma3} The operator $\Q:Z_T^0\to Z_T^0$ is compact.
\end{lemma}
\par
{\it Proof of Lemma  \ref{lemma3}}. Let $u=\L_0\Phi$, where $\Phi\in
Z_T^0$. By the semi-group property of backward SPDEs from Theorem
6.1 from Dokuchaev (2010), we obtain that
$u|_{t\in[0,s]}=\L_su(\cdot,s)$ for all $s\in (0,T]$. By Lemmas
\ref{lemma1} and \ref{Th3.1.1}, we have for $\tau\in\{t_1,...,t_m\}$
that \baaa \|\E\oo\G_iu(\cdot,\tau)\|^2_{W_2^1(D)}\le C_0
\|u(\cdot,\tau)\|^2_{Z^1_0}\le C_1\inf_{t\in[\tau,T]} \|
u(\cdot,t)\|^2_{Z_t^1}\le \frac{C_1}{T-\tau}
\int_\tau^T\|u(\cdot,t)\|_{Z_t^1}^2dt\\\le
\frac{C_2}{T-\tau}\|\Phi\|^2_{Z_T^0}\eaaa and \baaa
 &&\|\E\oo\G_0u\|^2_{W_2^1(D)}\le 
C_3\E\int_0^T\|u(\cdot,t)\|^2_{Z^1_T}dt\le C_4 \|\Phi\|_{Z_T^0}.
 \eaaa for constants $C_i>0$ which
do not  depend on $\Phi$. Hence the operator $\Q:Z_T^0\to W_2^1(D)$
is continuous. Since the embedding of $W_2^1(D)$ to $H^0$ and in
$Z_T^0$ is a compact operator, the proof of Lemma \ref{lemma3}
follows. $\Box$
\par
{\it Proof of  Theorem \ref{Th3}}. By the assumptions, the equation
$\Q \Phi=\Phi$ has the only solution $\Phi=0$ in $H^0$. By Lemma
\ref{lemma3} and by the Fredholm Theorem, the operator
$(I-\Q)^{-1}:H^0\to H^0$ is continuous. Then the proof of  Theorem
\ref{Th3} follows from representation (\ref{4.3}). $\Box$
\par
Let us introduce operators
   $$ \A^* v\defi\sum_{i,j=1}^n\frac{\p^2 }{\p x_i \p x_j}
\left(b_{ij}(x,t)v(x)\right)-\sum_{i=1}^n\frac{\p}{\p x_i }\left( \w
f_i(x,t)v(x)\right)+\w\lambda(x,t)v(x)$$ and \baaa B_i^*v\defi
 -\sum_{k=1}^n \frac{\p }{\p x_k }\,\big(\beta_{ik}(x,t,\o)\,v(x))+
 \oo \beta_i(x,t,\o)\,v(x),\qquad i=1,\ldots ,N. \label{AB*}\eaaa
Here $b_{ij}$, $x_i$, $\b_{ik}$ are the components of $b$,  $\b_i$,
and $x$.
\par
Let $\rho\in Z_{s}^0$, and let $p=p(x,t,\o)$ be the solution of  the
problem \baaa &&d_tp=\A^* p\, dt +
\sum_{i=1}^NB^*_ip\,dw_i(t), \quad t\ge s,\nonumber\\
&&p|_{t=s}=\rho,\quad\quad p(x,t,\o)|_{x\in \p D}=0.\label{p}\eaaa
By Theorem 3.4.8 from Rozovskii (1990), this boundary value problem
has an unique solution $p\in Y^1(s,T)$. Introduce an operator $\M_s:
Z_s^0\to Y^1(s,T)$  such that $p=\M_s\rho$, where $p\in Y^1(s,T)$ is
the solution of this boundary value problem.
\par
{\em Proof of  Theorem \ref{Th5}.}
By Theorem 3.1 from Dokuchaev (2005), problem (\ref{4.1})  has an unique solution $p\in {Y}^2$ for any
$\rho \in Z_s^1$, and \be \| p\|_{{Y^2(s,T)}}
 \le   C
\|\rho \|_{Z_s^1}, \label{3.1f} \ee where  $C>0$ does not depend on
$\rho$ (Dokuchaev (2005)). This $C$ depends on ${\cal P}$ and on the supremums of the
derivatives in Condition \ref{condB}.
\par
By Theorem 4.2 from Dokuchaev (2010), we have that
$\kappa p(\cdot,T)=\Q^*\rho$, i.e., \baa
(\rho,\Q\Phi)_{Z_0^0}=(\rho,\kappa v(\cdot,0))_{Z_0^0}=
(p(\cdot,T),\kappa v(\cdot,T))_{Z_T^0}=
(\kappa p(\cdot,T),\Phi)_{Z_T^0}\label{dual}\eaa for $v=\L_T\Phi$. (See
also Lemma 6.1 from Dokuchaev (1991) and related results in Zhou
(1992)). \par Suppose that there exists $\Phi\in Z_T^0$ such that
$\kappa v(\cdot,0)=v(\cdot,T)$ for $v=\L_T\Phi$, i.e.,
$v(\cdot,0)=\Q\Phi=\Phi$.
Let us show that $\Phi=0$ in this case.
\par
\def\QQ{{\bf Q}}
Since $\Q\Phi\in Z_0^0$, it follows that $\Phi\in H^0=Z_0^0$. Let $p=\M_0\rho$ and
$\oo p(x,t,0)=\E p(x,t,\o)$ (meaning the projection from $Z_T^0$ on
$H^0=Z_0^0$). Introduce an operator $\QQ: H^0\to H^0$  such that
$\kappa\oo p(\cdot,T)=\QQ\rho$. By (\ref{dual}), the properties of $\Phi$ lead to the equality \baa
(\rho-\kappa p(\cdot,T),\Phi(\cdot,T))_{Z_T^0}=(\rho-\kappa \oo
p(\cdot,T),\Phi(\cdot,T))_{H_0}=0 \quad\forall \rho\in
H^0.\label{pPhi}\eaa It suffices to show that the set $\{\rho-\kappa \oo
p(\cdot,T)\}_{\rho\in H^0}$ is dense in $H^0$. For this, it suffices
to show that the equation $\rho-\QQ\rho=z$ is solvable in $H^0$ for
any  $z\in H^0$.
\par
Let us show that the operator $\QQ:H^0\to H^0$ is compact. Let $p$
be the solution of (\ref{p}). This means that $\kappa \E p(\cdot,T)=\QQ
\rho$. By Lemma \ref{Th3.1.1}, it follows that \be \label{2f}
\|p(\cdot,\tau)\|_{Z_\tau^1}\le C\|p(\cdot,s)\|_{Z_s^1},\quad
\tau\in[s,T], \ee where $C_*>0$ is a constant that does not depend
on $p$, $s$, and $\tau$.
\par
We have that $p|_{t\in[s,T]}=\M_s p(\cdot,s)$ for all $s\in
[0,T]$, and, for $\tau>0$, \baaa \|\oo p(\cdot,T)\|^2_{W_2^1(D)}&\le
& C_0\|p(\cdot,T)\|^2_{Z^1_T}\le C_1\inf_{t\in[0,T]}\|
p(\cdot,t)\|^2_{Z_t^1}\\ &\le& \frac{C_1}{T} \int_0^T\|
p(\cdot,t)\|^2_{Z_t^1}dt  \le \frac{C_2}{T}\|p\|^2_{X^1} \le
\frac{C_3}{T}\|\Phi\|_{H^0} \eaaa for constants $C_i>0$ that do not
depend on $\Phi$. Hence the operator $\QQ:H^0\to H^1$ is continuous.
The embedding of $H^1$ into $H^0$ is a compact operator (see, e.g.,
Theorem 7.3 from Ladyzhenskaia (1985), Chapter I).
\par
Similar to the proof of Theorem 3.4 from Dokuchaev (2008), it can be
shown that if \baa\kappa \oo p(\cdot,T)=\kappa \E
p(\cdot,T)=\QQ\rho=p(\cdot,0)\label{pp}\eaa for some $\rho\in H^0$
then $\rho=0$; this proof will be omitted here.
 We had proved also that the operator $\QQ$ is compact. By
the Fredholm Theorem, it follows that the equation $\rho-\QQ\rho=z$
is solvable in $H^0$ for any $z\in H^0$. By (\ref{pPhi}), it follows
that $\Phi=0$.
 Therefore, the condition $\kappa u(\cdot,0)=u(\cdot,T)$
fails to be satisfied for $u\neq 0$, $\xi=0$, and $\varphi=0$. Thus, $u=0$ is the
unique solution of problem (\ref{parab1})-(\ref{parab2}) for $\xi=
0$ and $\varphi=0$. Then the proof of Theorem \ref{Th5} follows
from Theorem \ref{Th3}. $\Box$
\par
Without a loss of generality, we assume that there exist functions
${\ww\b_i: Q\times \O \to \R^n}$, $i=1,\ldots, M$, such that $$
2b(x,t,\o)=\sum_{i=1}^N\b_i(x,t,\o)\,\b_i(x,t,\o)^\top
+\sum_{j=1}^M\,\ww\b_j(x,t,\o)\,\ww\b_j(x,t,\o)^\top, $$ and $\ww
\b_i$ has the similar properties as $\b_i$. (Note that, by Condition
\ref{cond3.1.A}, $2b\ge \sum_{i=1}^N\b_i\b_i^\top$).
\par
 Let
$\ww w(t)=(\ww w_1(t),\ldots, \ww w_M(t))$ be a new Wiener process
independent on $w(t)$. Let $a\in L_2(\O,\F,\P;\R^n)$ be a vector
such that $a\in D$. We assume also that $a$ is independent from
$(w(t)-w(t_1),\w w(t)-\w w(t_1))$ for all $t>t_1>s$. Let $s\in[0,T)$
be given. Consider the following Ito equation
\begin{eqnarray}
\label{yxs} &&dy(t) =
\w f(y(t),t)\,dt+\sum_{i=1}^N\b_i(y(t),t)\,dw_i(t) +\sum_{j=1}^M\ww
\b_j(y(t),t)\,d \ww w_j(t),
\nonumber\\ [-6pt] &&y(s)=x.
\end{eqnarray}
\par
Let  $y(t)=y^{x,s}(t)$ be the solution of (\ref{yxs}), and let
$\tau^{x,s}\defi\inf\{t\ge s:\ y^{a,s}(t)\notin D\}$.
 For $t\ge s$,
set
\baa
\g^{x,s}(t) \defi\exp\left(-\int_s^t
\w\lambda(y^{x,s}(t),t)\,dt\right).\label{gamma}
\eaa
\begin{lemma}\label{ThJust} Let Conditions \ref{condB}-\ref{condTh5} hold.
  Let $\Phi\in L_{\infty}(\O,\F_T,\P,C_0(\oo D))$,     and let   $u\in Y^1$ be solution of
(\ref{4.1})  with $\varphi=0$. Then the process $\g^{x,s}(t\land \tau)u(y^{x,s}(t\land \tau),t\land
\tau)$ is a martingale.
\end{lemma}
\par
 {\it Proof of Lemma \ref{ThJust}}. For the case where
$u\in\X_c^2$ and $\chi_j\in\X_c^1$, this lemma follows  from the
proof of Lemma 4.1 from Dokuchaev (2011) (see Remark \ref{remDu}).
\par
Let us consider the general case.
 Let $\rho\in Z_s^0$ be such that $\rho\ge 0$ a.e. and $\int_{D}\rho(x)dx=1$ a.s.
 Let $a\in L_2(\O,\F,\P;\R^n)$ be
  such that $a\in D$ a.s. and it
has the conditional probability  density function $\rho$ given
$\F_s$. We assume that $a$ is independent from  $(w(t_1)-w(t_0),\w
w(t_1)-\w w(t_0)\}$, $s<t_0<t_1$. Let $p=\M_s\rho$, and let
  $y^{a,s}(t)$ be
the solution of  Ito equation (\ref{yxs}) with the initial condition
$y(s)=a$.
\par
To prove the theorem, it suffices to show that \baaa
\g(t\land\tau^{x,s})u(y^{x,s}(t\land\tau^{x,s}),t\land\tau^{x,s})=\E_t\g(T\land\tau^{x,s})
u(y^{x,s}(T\land\tau^{x,s}),T\land\tau^{x,s})\quad \hbox{a.s.}\eaaa for any $t$. For this, it suffices to prove that \baa
&&\E\int_D\rho(x)
\g(t\land\tau^{x,s})u(y^{x,s}(t\land\tau^{x,s}),t\land\tau^{x,s})\nonumber\\&&=\E\int_D\rho(x)\E_t
\g(T\land\tau^{x,s})u(y^{x,s}(T\land\tau^{x,s}),T\land\tau^{x,s})\label{ident}\eaa
 for any $\rho\in
Z_s^0$ such as described above.
\par
By Theorem 6.1 from Dokuchaev (2011) and Remark \ref{remDu}, we have that \baaa
\int_Dp(x,t)u(x,t)dx=\E_t
\g^{a,s}(t\land\tau^{a,s})u(y^{a,s}(t\land\tau^{a,s}),t\land\tau^{a,s})\eaaa
and \baaa \int_Dp(x,T)u(x,T)dx=\E_T\,\g(T\land\tau^{a,s})
u(y^{a,s}(T\land\tau^{a,s}),T\land\tau^{a,s}).\eaaa By the duality
established in Theorem 3.3 from Dokuchaev (2011) and Remark \ref{remDu}, it follows that
\baaa\E\int_Dp(x,t)u(x,t)dx=\E\int_Dp(x,T)u(x,T)dx.\eaaa
This means that $\E(\E_tq(a,s,t))=\E(\E_Tq(a,s,T))$,
where $$q(a,s,t)=\g^{a,s}(t\land\tau^{a,s})u(y^{a,s}(t\land\tau^{a,s}),t\land\tau^{a,s}).$$
Hence \baa
\E(\E_tq(a,s,t))=\E(\E_tq(a,s,T)).\label{eqeq}\eaa
\par
Without loss of generality, we shall assume that $a$ is a random
vector on the probability space $(\ww\O,\ww\F,\ww\P)$, where
$\ww\O=\O\times\O'$, where $\O'=D$,     $\ww\F=\overline{\F_s\otimes
{\cal B}_D}$, where ${\cal B}_D$ is the set of Borel subsets of $D$,
and \baaa \ww\P(S_1\times S_2)=\int_{S_1}\P(d\o)\P'(\o,S_2),\quad
\P'(\o,S_2)=\int_{S_2}\rho(x,\o)dx, \eaaa for $S_1\in\F$ and $S_2\in
{\cal B}_D$. The symbol $\ww\E$ denotes the expectation in
$(\ww\O,\ww\F,\ww\P)$. We suppose that $\ww\o=(\o,\o')$,
$\ww\O=\{\oo\o\}$, and $a(\ww w)=\o'$.
\par
We have that
\baaa
\E(\E_t q(a,s,t))=\E\int_{\ww\O}\ww\P(d\o|\F_t)q(\o',s,t,\o)=\E\int_Dd\o'\rho(\o')\int_{\O}\P(d\o|\F_t)q(\o',s,t,\o)\\=\E\int_D\E_t\rho(\o') q(\o',s,t,\o)d\o'=\E\int_D\rho(\o')q(\o',s,t,\o)d\o'\\
=\E\int_D\rho(x) \g^{x,s}(t\land\tau^{x,s})u(y^{x,s}(t\land\tau^{x,s}),t\land\tau^{x,s})dx\eaaa
and
\baaa
\E(\E_tq(a,s,T))=\E\int_{\ww\O}\ww\P(d\o|\F_t)q(\o',s,T,\o)=\E\int_Dd\o'\rho(\o')\int_{\O}\P(d\o|\F_t)q(\o',s,T,\o)
\\=\E\int_D\rho(\o')\E_t q(\o',s,T,\o)d\o'
=\E\int_D\rho(x)
\g^{x,s}(T\land\tau^{x,s})\E_t u(y^{x,s}(T\land\tau^{x,s}),T\land\tau^{x,s})dx.
\eaaa  Since the choice of $\a$ and $\rho$ is arbitrarily, it
follows from (\ref{eqeq})  that (\ref{ident}) holds. This completes
the proof of Lemma \ref{ThJust}. $\Box$
\par
\def\ttau{T\land\tau}
{\em Proof of Theorem \ref{ThF}}. Without a loss of generality, we
assume that $\P$ is a  martingale probability measure, i.e., $S(t)$
is a martingale and $dS(t)=\s(t)S(t)dw(t)+\w\s(t)S(t)d\w w(t)$. Let
$\zeta=u(S(T\land\tau),T\land\tau)$.

 We have that $u\in Y_1$. Hence $u(\cdot,T)=Z_T^1$. Similarly to the proof of  Lemma  \ref{lemma3}, we obtain that
$u(\cdot,0)\in Z_0^1$, and  $\|u(\cdot,0)\|_{Z_0^1}\le
C\|u(\cdot,T)\|_{Z_T^0}$, where $C>0$ does not depend on
$u(\cdot,T)$. Since the embedding of $H^1$ to $C_0(\oo D)$ is
continuous for $n=1$, we obtain that $u(\cdot,0)\in
L_2(\O,\F,\P,C_0(\oo D))$. By (\ref{XT}) and by the assumptions on
$\xi$, we obtain that $u(\cdot,T)\in L_2(\O,\F,\P,C(\oo D))$. Hence
$\E \zeta^2\le \E (\sup_{x}|u(x,T)|^2)<+\infty$.
\par
Since the market is complete, there exists admissible $\g(t)$ such
that \baaa \zeta=u(S(\ttau),\ttau)=\E
u(S(\ttau),\ttau)+\int_0^T\g(t)dS(t). \eaaa   Therefore,
$\E\{u(S(\ttau),\ttau)|\F_t\}$ is the wealth for the self-financing
strategy such that replicates $\zeta$.

By Lemma \ref{ThJust}, it follows that
$\E\{u(S(\ttau),\ttau)|\F_t\}=u(S(t\land \tau),t\land\tau).$ Hence
(\ref{Xmar}) is the wealth for the self-financing strategy that
replicates $\zeta$.
\par
Let $S_1(t)=S(t)/S(0)$. Further, if $\tau>T$ and $S(T)y=S(0)x$,
then, by the definitions, \baaa && X(T,y)-X(0,x)=u(S_1(T)y,T)
-u(S_1(0)x,t) = u(x,T) -u(x,0)=\xi(x). \eaaa Hence condition
(\ref{XT}) is satisfied. This completes the proof of Theorem
\ref{ThF}. $\Box$

\index{\begin{theorem}\label{ThF} $X(T)=X(0)+\int_0^\tau \varphi(S(t),t)dt$
The investment problem  (\ref{XL})-(\ref{XT}) has a solution with the terminal
 wealth  $X(T)=u(S(T\land\tau),T\land\tau)+\int_0^{T\land\tau}\varphi(S(T),t)$.
\end{theorem}}

\subsection*{Acknowledgment} This work  was
supported by ARC grant of Australia DP120100928 to the
author.
\section*{References} $\hphantom{XX}$
Al\'os, E., Le\'on, J.A., Nualart, D. (1999).
 Stochastic heat equation with random coefficients
 {\it
Probability Theory and Related Fields} {\bf 115} (1), 41--94.
\par
Bally, V., Gyongy, I., Pardoux, E. (1994). White noise driven
parabolic SPDEs with measurable drift. {\it Journal of Functional
Analysis} {\bf 120}, 484--510.
\par Caraballo, T.,  Kloeden, P.E.,
Schmalfuss, B. (2004). Exponentially stable stationary solutions for
stochastic evolution equations and their perturbation, Appl. Math.
Optim. {\bf  50}, 183--207.
\par
 Chojnowska-Michalik, A. (1987). On processes of Ornstein-Uhlenbeck type in
Hilbert space, Stochastics 21, 251--286.

\par
Chojnowska-Michalik, A. (1990). Periodic distributions for linear
equations with general additive noise, Bull. Pol. Acad. Sci. Math.
38 (1–12) 23--33.

\par
Chojnowska-Michalik, A., and Goldys, B. (1995). {Existence,
uniqueness and invariant measures for stochastic semilinear
equations in Hilbert spaces},  {\it Probability Theory and Related
Fields},  {\bf 102}, No. 3, 331--356.
\par
Da Prato, G., and Tubaro, L. (1996). { Fully nonlinear stochastic
partial differential equations}, {\it SIAM Journal on Mathematical
Analysis} {\bf 27}, No. 1, 40--55.
\par
Dokuchaev, N.G. (1992). { Boundary value problems for functionals
of
 Ito processes,} {\it Theory of Probability and its Applications}
 {\bf 36} (3), 459-476.
\par
 Dokuchaev, N.G. (1994). Parabolic equations without the Cauchy
boundary condition and problems on control over diffusion processes.
I. {\em Differential Equations} {\bf 30}, No. 10, 1606-1617;
translation from Differ. Uravn. 30, No.10, 1738-1749.
 \par Dokuchaev, N.G. (2004).
Estimates for distances between first exit times via parabolic
equations in unbounded cylinders. {\it Probability Theory and
Related Fields}, {\bf 129}, 290 - 314.
\par
Dokuchaev, N.G. (2005).  Parabolic Ito equations and second
fundamental inequality.  {\it Stochastics} {\bf 77} (2005), iss. 4.,
pp. 349-370.
\par
 Dokuchaev N. (2008) Parabolic Ito equations with mixed in time
conditions.
{\it Stochastic Analysis and Applications} {\bf 26}, Iss. 3, 562--576. 
\par
Dokuchaev, N. (2010). Duality and semi-group property for backward
parabolic Ito equations. {\em Random Operators and Stochastic
Equations. } {\bf 18}, 51-72.
\par
Dokuchaev, N. (2011). Representation of functionals  of Ito
processes in bounded domains. {\em Stochastics} {\bf 83}, No. 1,
45--66.
\par
 Dokuchaev, N. (2012a).
Backward parabolic Ito equations and second fundamental inequality.
{\em Random Operators and Stochastic Equations} {\bf 20}, iss. 1,
69-102.
\par
Dokuchaev, N, (2012b). Backward SPDEs with non-local in time and
space boundary conditions. Working paper, arXiv:1211.1460 (submitted).
\par
 Dokuchaev, N. (2012c). On almost surely
periodic and almost periodic solutions of backward SPDEs.  Working paper, arXiv: 1208.5538 (submitted).
\par
Du K., and Tang, S. (2012). Strong solution of backward stochastic
partial differential equations in $C^2$ domains. {\em  Probability
Theory and Related Fields)} {\bf 154}, 255--285.
\par
Duan J.,  Lu K., Schmalfuss B. (2003). Invariant manifolds for
stochastic partial differential equations.{\em Ann. Probab.} {\bf
31}
 2109–2135.
\par
Feng C., Zhao H. (2012). Random periodic solutions of SPDEs via
integral equations and Wiener-Sobolev  compact embedding. {\em
Journal of Functional Analysis} {\bf 262}, 4377--4422.
\par
Gy\"ongy, I. (1998). Existence and uniqueness results for semilinear
stochastic partial differential equations. {\it Stochastic Processes
and their Applications} {\bf 73} (2), 271-299.
\par
Kl\"unger, M. (2001). Periodicity and Sharkovsky's theorem for
random dynamical systems, {\em Stochastic and Dynamics} {\bf 1},
iss.3, 299--338.
\par Krylov, N. V. (1999). An
analytic approach to SPDEs. Stochastic partial differential
equations: six perspectives, 185--242, Mathematical Surveys and
Monographs, {\bf 64}, AMS., Providence, RI, pp.185-242.
\par
Ladyzhenskaia, O.A. (1985). {\it The Boundary Value Problems of
Mathematical Physics}. New York: Springer-Verlag.
\par
Liu, Y., Zhao, H.Z (2009). Representation of pathwise stationary
solutions of stochastic Burgers equations, {\em Stochactics and
Dynamics} {\bf  9} (4), 613--634.
\par
 Ma, J., Yong, J. (1997). Adapted solution of a class
 of degenerate backward stochastic partial differential equations, with applications.
{\em Stochastic Processes and Their Applications.} Vol. 70, pp.
59--84.
\par
Maslowski, B. (1995). { Stability of semilinear equations with
boundary and pointwise noise}, {\it Annali della Scuola Normale
Superiore di Pisa - Classe di Scienze} (4), {\bf 22}, No. 1,
55--93.
\par
Mattingly. J. (1999). Ergodicity of 2D Navier–Stokes equations with random forcing and large viscosity. {\em Comm. Math.
Phys.} 206 (2),  273–288.
\par
Mohammed S.-E.A., Zhang T.,  Zhao H.Z. (2008). The stable manifold
theorem for semilinear stochastic evolution equations and stochastic
partial differential equations. {\em Mem. Amer. Math. Soc.} 196
(917),  1–105.
\par
Pardoux, E. (1993). Bulletin des Sciences Mathematiques, 2e Serie,
 {\bf 117}, 29-47.
\par
 Revuz, D., and Yor, M. (1999). {\it Continuous Martingales
and Brownian Motion}. Springer-Verlag: New York.
\par
Rodkina, A.E. (1992). On solutions of stochastic equations with
almost surely periodic
    trajectories.  {\em Differ. Uravn}. 28, No.3, 534--536 (in Russian).
\par
Rozovskii, B.L. (1990). {\it Stochastic Evolution Systems; Linear
Theory and Applications to Non-Linear Filtering.} Kluwer Academic
Publishers. Dordrecht-Boston-London.
\par
Sinai, Ya. (1996). Burgers system driven by a periodic stochastic
flows, in: Ito's Stochastic Calculus and Probability Theory,
Springer, Tokyo, 1996, pp. 347–353.
  \par
Walsh, J.B. (1986). An introduction to stochastic partial
differential equations, Lecture Notes in Mathematics {\bf 1180},
Springer Verlag.
\par
Yong, J., and Zhou, X.Y. (1999). { Stochastic controls: Hamiltonian
systems and HJB equations}. New York: Springer-Verlag.
\par
 Zhou, X.Y. (1992). { A duality analysis on stochastic partial
differential equations}, {\it Journal of Functional Analysis} {\bf
103}, No. 2, 275--293.
\end{document}